\documentclass[10pt]{article}

\usepackage{epsfig}
\usepackage{amsmath}
\usepackage{amssymb}
\usepackage{graphicx}
\usepackage{color}

\usepackage[latin1]{inputenc}

\newcommand{\R}{\mathbb{R}}

\def\bbI{{\rm 1\!\!{\rm I}}}
\def\phi{\varphi}

\def\bbr{{\mathbb R}}

\def\bz{{\bf Z}}

\def\bz{{\mathbb Z}} %Keistas Z apribrezimas

\def\bk{{\bf k}}
\def\bt{{\bf t}}
\def\bs{{\bf s}}

\def\bX{{\bf X}}
\def\bY{{\bf Y}}

\def\Bl1{{\bf 1}}
\def\B2{{\bf 2}}
\def\B0{{\bf 0}}

\def\a{\alpha}
\def\b{\beta}
\def\d{\delta}

\def\e{\varepsilon}
\def\g{\gamma}
\def\G{\Gamma}
\def\l{\lambda}

\def\s{\sigma}
\def\t{\tau}

\def\=A8{\"o}
\def\lygpas{\stackrel{ D}{\ =\ }}

\newcommand{\beq}{\begin{equation}}
\newcommand{\eeq}{\end{equation}}
\newcommand\beqn{\begin{displaymath}}  % no number
\newcommand\eeqn{\end{displaymath}}

\newcommand{\halmos}{\vspace{3mm} \hfill \mbox{$\Box$}\\[2mm]}
\newtheorem{teo}{Theorem}

\newtheorem{cor}[teo]{Corollary}
\newtheorem{prop}[teo]{Proposition}
\newtheorem{definition}[teo]{Definition}

\begin{document}

\title{On $\a$-covariance, long, short and negative memories for
sequences of random variables  with infinite variance}
\author{ Vygantas Paulauskas}
\date{July 17, 2013 \\  \small
\vskip.2cm Vilnius University }
 \maketitle

\begin{abstract}
We consider a measure of dependence for  symmetric $\a$-stable
random vectors, which was introduced by the author in 1976. We
demonstrate that this measure of dependence can be extended for much
more broad class of random vectors (up to regularly varying vectors
in separable Banach spaces). This measure is applied for linear
random processes and fields with heavy-tailed  innovations, for some
stable processes, and these applications show that this dependence
measure, named as $\a$-covariance is a good substitute for the usual
covariance.

Also we discuss a problem of  defining  long, short, and negative
memories for stationary processes and fields with infinite
variances.

\end{abstract}
\vfill
\eject
\section{Introduction }
%\subsection{Introduction}

%\subsection{Results}
  The importance of the notion of independence in
 probability theory is well-known, sometimes it is even stressed that this is the
 main feature which distinguishes the probability from the general
 measure theory. Therefore, the notion of dependence in probability
 theory  as being
 in some sense opposite to independence, is also important,
 moreover, it is much more complicated for the following reason.
 Let $X$ and $Y$ be two random variables (or a random vector $(X, Y)$)
  defined on some probability space. They are independent if their
  joint distribution is a product of marginal distributions. But if
  they are not independent, we would like to know what kind of
  dependence between $X$ and $Y$ we are facing with, how strong
  this dependence is. This means that in the case of dependence we
  want to measure this dependence, and we would like  to have the
  property that for independent random variables this measure
   would be zero, while for the "strongest" dependence it would be one.
  Here "measure" is used not in the measure-set theoretical meaning,
  here it stands for some  function (even not  necessarily non-negative,
  as in the case of a correlation coefficient, which can be both positive and negative)
   on the set of all bivariate distributions.
  Again, what is the "strongest" dependence it is not quite
  clear, one possible candidate for such dependence can be the case
  where one random variable is a function of another, i.e., $Y=f(X)$
  where $f$ is one-to-one function from a support of $X$ to a
  support of $Y.$ In probability theory and mathematical statistics
   a lot of measures or concepts of dependence are
  introduced, among them classical Pearson correlation coefficient,
  Kendall's $\tau$, Spearman's $\rho$, more recent functional or physical dependence measure,
   and many others. We refer for
   recent survey papers \cite{Bradley1} and \cite{Wu}  with  big lists of references for
  measures of dependence.

In the case of random variables having finite second moments one of
the most popular measures of dependence is correlation coefficient,
defined by the following formula
$$
Corr (X,Y)=\frac{E(X-EX)(Y-EY)}{\sigma(X) \sigma(Y)},
$$
where $\sigma ^2(X)=E(X-EX)^2.$ For independent random variables
correlation coefficient is zero, always $-1\le Corr (X,Y)\le 1$ and
if $Y=aX+b$ for some real numbers $a, b$, then $|Corr (X,Y)|=1.$
Unfortunately, equality $Corr (X,Y)=0$ holds not only for
independent $X$ and $Y$, it can be even for random vector $(X,Y)$
concentrated on some curve, the case, which we would like to
attribute as the "strongest" dependence. Thus, it is possible to say
that correlation coefficient measures linear dependence and random
variables satisfying $Corr (X,Y)=0$ are called uncorrelated. This
notion is very important in the so-called $L_2$-theory of random
variables (uncorrelated means orthogonal in this theory). These
properties and the simplicity of the notion explain why measures of
dependence based on correlation and covariation are so popular and
are used in many areas of probability and statistics, in particular,
in time series analysis.  One of the main ways to define memory
properties (long, short, and negative memories) for a covariance
stationary mean zero process $Y(t)$ is to use the decay and some
other properties of the covariance function, see the precise
Definition \ref{def1}. In its turn, these properties are important
when considering functional limit theorems for these processes. One
class of stationary processes, for which these memory properties are
studied most deeply are linear processes. Let $\e_i, \ i\in \bz,$ be
independent, identically distributed (i.i.d.) random variables with
finite second moment (without loss of generality we may assume
$E\e_1=0, E\e_1^2=1$) and let $a_k, k\ge 0$ be a sequence of real
numbers satisfying $\sum_{k=0}^\infty a_k^2<\infty$ (this sequence
sometimes is called a filter of a linear process, while random
variables $\e_i, \ i\in \bz,$ are called innovations). Then a linear
process
$$
X_t=\sum_{k=0}^\infty a_k \e_{t-k}, \ t\in \bz,
$$
is a stationary sequence, and the dependence is  reflected in the
covariance function
$$
\gamma(k)= EX_0X_k, \ k\in \bz.
$$
Since there is a simple expression of $\gamma(k)$ via coefficients
of the filter ($\gamma(k)=\sum_{j=0}^\infty a_ja_{j+k}$), properties
of the filter define  memory properties of the linear process under
consideration. In Section 4 we discuss this question in detail and
argue that memory properties can not be defined only by the decay of
covariance function.

The situation is quite different if random variables $X$ and $Y$
have infinite variance, one can say that in this case there is no
good substitute for correlation coefficient. Of course, in
statistics there are above mentioned Kendall's $\tau$ or Spearman's
$\rho$, which are based on rank statistics and, therefore, do not
require any moments of random variables under consideration. But
such statistics are not convenient for investigation dependence in
more theoretical problems, for example, they are not convenient to
measure dependence in linear processes with innovations without
second moment. On the other hand, during last decades the role of
the so-called heavy tailed distributions had increased both in
theoretical and applied probability, therefore the problem of
measuring the dependence between random variables having infinite
variance remains an important (and difficult) problem.

The main aim of this paper is to revive interest to one measure of
dependence which was introduced by the author more than 30 years ago
in \cite{Paul3} for a specific, but rather important, class of
 symmetric $\a$ -stable
($S\a S$) distributions. The importance of this class of
distributions can be explained by the fact that trying to build a
model involving random variables with infinite variance, as a first
step,  one takes stable random variables, or random variables in the
domain of their attraction. Although at the beginning we shall deal
mainly with bivariate random vectors, we shall define general
$d$-dimensional $S\a S$  random vectors. Let $S_d=\{x\in \R^d:
||x||=1\}$ be the unit sphere in $\R^d$, here $||x||$ stands for the
Euclidean norm in $\R^d$. Random vector $\bX =(X_1, \dots X_d)$ is
$S\a S$ with parameter $0<\a<2$ if there exists a unique symmetric
finite measure $\G$ on  $S_d$ such that the characteristic function
(ch.f.) of $\bX$ is given by formula
\begin{equation}\label{stabdef}
E\exp \left \{i(\bt, \bX) \right \}=\exp \left \{-\int_{S_d}|(\bt,
\bs)|^\a \G (d\bs) \right \}.
\end{equation}
$\G$ is called the spectral measure of the $S\a S$ random vector
$\bX$. The Gaussian case $\a=2$ is excluded from this definition,
since in the Gaussian case there is no uniqueness of the spectral
measure $\G$: many different measures $\G$ will give the same ch.f..
About forty years ago (for historical details we refer to monograph
\cite{Samorod}) two measures of dependence between coordinates of a
bivariate $S\a S$ random vector $\bX =(X_1, X_2)$ with spectral
measure $\G$ were introduced. The first one, called covariation of
$X_1$ on $X_2$ and denoted  by $[X_1, X_2]_\a$, is defined for
$\a>1$ as follows:
$$
[X_1, X_2]_\a =\int_{S_2}s_1s_2^{\langle \a-1\rangle} \G (d\bs),
$$
where $a^{\langle p\rangle}= |a|^p \ {\rm sign}\ a.$ Although
 in the case
$\a=2$ this quantity is equal to the half of covariance between
$X_1$ and $X_2$, which is symmetric,  it is  not symmetric in its
arguments, i.e., in general (for $1<\a<2$)
$$
[X_1, X_2]_\a \ne [X_2, X_1]_\a
$$
All properties of covariation, including equivalent definition, are
given in Chapter 2.7 of \cite{Samorod}. Main shortcoming of this
measure of dependence, apart of the just mentioned non-symmetricity,
is that it is not defined for $\a<1$ (for $\a=1$ it is possible to
define the covariation, see Exercise 2.22 in \cite{Samorod}). Here
it is appropriate to mention that recently in \cite{Garel1}
symmetric  covariation coefficient was introduced.

Another measure of dependence for $S\a S$ random vectors is the
codifference, defined by formula
$$
\t (X_1, X_2) =\int_{S_2}\left (|s_1|^\a +|s_2|^\a -|s_1-s_2|^\a
\right ) \G (d\bs).
$$
One can note that the codifference can be defined for a general
bivariate random vector $\bY =(Y_1, Y_2)$  by means of the following
formula
\begin{equation}\label{codifdefgen}
\t (Y_1, Y_2) =\ln f_{\bY}(1,-1) -\ln f_{\bY}(1,0) -\ln
f_{\bY}(0,-1),
\end{equation}
where
$$
f_{\bY}(t,s)= E\exp\{i(tY_1+s Y_2)\}.
$$
 The codifference for $S\a S$ random vectors has better properties: it
is symmetric function ($\t (X_1, X_2)=\t (X_2, X_1)$), is defined
for all $0<\a\le 2$, in the Gaussian case coincides with the
covariance. Other properties one can find in Chapter 2.10 of
\cite{Samorod}. Here it is appropriate to mention the recent papers
\cite{Szekely1} and \cite{Szekely2}, where the distance covariance
and Brownian distance covariance where introduced, and these
measures of dependence are based on the same idea as the
codifference: ch.f. of a vector $\bY$ with independent coordinates
is a product of marginal ch.f. of components.

Earlier than the codifference and about the same time as covariation
were introduced, the author in \cite{Paul3} had proposed one more
measure of dependence for $S\a S$ random vectors. Let $X=(X_1, X_2)$
be a $S\a S$ random vector with the spectral measure $\G$, and let
$Y=(Y_1, Y_2)$ be a random vector on $S_2$ with the distribution
${\tilde \G}(A)= (\G(S_2))^{-1}\G(A)$. Then the generalized
association parameter (g.a.p.) of the random vector $X$ is defined
as usual correlation coefficient for the random vector $Y$:
\begin{eqnarray} \label{gapdef}
{\tilde \rho}={\tilde \rho} (X_1, X_2) &=&
\frac{EY_1Y_2}{\sqrt{EY_1^2EY_2^2}}=\frac{\int_{S_2}s_1s_2{\tilde
\G}(ds)}{\left (\int_{S_2}s_1^2{\tilde \G}(ds)\int_{S_2}s_2^2{\tilde
\G}(ds) \right )^{\frac{1}{2}}} \\ \nonumber
&=&\frac{\int_{S_2}s_1s_2{\G}(ds)}{\left
(\int_{S_2}s_1^2{\G}(ds)\int_{S_2}s_2^2{\G}(ds) \right
)^{\frac{1}{2}}}.
\end{eqnarray}
Also we shall use the following analog of the covariance between
$X_1$ and $X_2$
\begin{equation} \label{acovdef}
{\rho}={\rho} (X_1, X_2) = \int_{S_2}s_1s_2{\G}(ds),
\end{equation}
and we shall call it as $\a$-covariance of a $S\a S$ random vector
$(X_1, X_2).$ Strictly speaking we should use the normalized measure
${\tilde \G}$ instead of $\G$ in the definition of $\a$-covariance
(then it would be possible to say that ${\rho} (X_1, X_2)$ is simply
covariance between random variables $Y_1$ and $Y_2$), but since $\G$
is finite measure, ${\rho} (X_1, X_2)= \G(S_2)EY_1Y_2.$

Here some remarks about the terminology is appropriate. The term
"generalized association parameter", clearly,  is  very
unsuccessful, it was introduced before the notions and terms
"covariation of $X_1$ on $X_2$" and "codifference" were invented, at
the time when the term "association" was popular (see, for example
papers \cite{Esary1}, \cite{Esary2}). Now we suggest to call this
parameter as $\a$-correlation coefficient for a $S\a S$ random
vector, and from now, in what follows we shall write
$\a$-correlation coefficient ($\a$-cc) instead of g.a.p. Thus we
have three notions (or quantities)- covariation, codifference, and
$\a$-covariance, which all become usual covariance in the case
$\a=2$.

 The following proposition was proved in \cite{Paul3}
\begin{prop}{\label{gapprop}}{\cite{Paul3}} The introduced $\a$-cc of
a random vector $\bX=(X_1, X_2)$ with ch.f. (\ref{stabdef}) (with
$d=2$) has the following properties:
\begin{description}
\item {(i)} $|{\tilde \rho}|\le 1$, and if the coordinates of $\bX$
are independent then ${\tilde \rho}=0$;
\item {(ii)} if $|{\tilde \rho}| = 1$, then the distribution of $\bX$
is concentrated on a line, i.e., coordinates $X_1$ and $X_2$ are
linearly dependent;
\item {(iii)} if $\a=2$, ${\tilde \rho}$ coincides with a
correlation coefficient of a Gaussian random vector with
characteristic function (\ref{stabdef});
\item {(iv)} ${\tilde \rho}$ is independent of $\a$ and depends only
on the spectral measure $\G$.
\end{description}
\end{prop}
Also in \cite{Paul3} it was shown that if a random vector $\bX$ is
sub-Gaussian with ch.f.
\begin{equation}\label{subgaus}
\exp \left \{-(\s_1^2t_1^2 +2r\s_1\s_2t_1t_2+\s_2^2t_2^2)^{\a/2}
\right \},
\end{equation}
where $\s_1^2, \ \s_2^2,$ and $r$ are variances and correlation
coefficient, respectively, of underlying Gaussian vector, then the
$\a$-cc ${\tilde \rho}=r.$

This notion is easily generalized to $d$-dimensional $S\a S$ random
vectors (see Proposition 3 in \cite{Paul3}) by defining the
$\a$-correlation matrix ${\tilde \Lambda}_\G$  and $\a$-covariance
matrix $ \Lambda_\G$ as usual correlation and covariance matrices,
respectively, of a random vector on $S_d$ with a distribution
${\tilde \G}(A)= (\G(S_d))^{-1}\G(A)$.

Despite of the simplicity of definition and the fact that $\a$-cc of
$S\a S$ random vectors satisfies main requirements for measures of
dependence, it was almost not used. Only recently the interest to
this measure of dependence was revived in \cite{Garel} and
\cite{Garel1}, where the so-called symmetric covariation coefficient
was introduced, it was compared with $\a$-cc, and estimations of
$\a$-cc and this new symmetric covariation  coefficient were
proposed. The main goal of the present paper is to demonstrate that
in the case of random vectors without variance these simple notions
of $\a$-covariance and $\a$-cc are quite natural substitutes for
covariance and correlation . We shall show that this measure of
dependence can be extended to more wide class of heavy-tailed random
vectors, is  very suitable when considering stable processes, in
particular, linear random processes and fields with heavy-tailed
innovations.

Also we discuss the memory properties for stationary sequences
without finite variance. We propose a unified approach to define
memory properties both for processes and for fields, based on the
rate of the growth of partial sums, formed by these processes or
fields. The exponent $1/\a$ (which is mentioned as boundary between
short and long memory in several papers) characterizes the growth of
partial sums of sequences with no memory at all, that is,  i.i.d.
random variables, and this value is attributed to short memory and
serves as a boundary between negative memory (when the exponent of
the growth is smaller than $1/\a$) and long memory (exponent is
bigger than $1/\a$).

By simple examples of linear processes and fields we show that this
approach is more natural even in the case of finite variance. We
want to stress that it is important to separate notions of memory
and dependence, that is, to separate long and short-range dependence
from memory properties, therefore it would be logical to call these
properties negative, zero and positive memories, and even we suggest
to introduce strongly negative memory (the case where the volatility
of partial sums stays bounded and do not grow with number of
summands in partial sums).

At the same time we agree with the attitude propagated by G.
Samorodnitsky in his several papers (see, for example,
\cite{Samorodnitsky} and \cite{Samorodnitsky1}), that memory
phenomenon is a complicated one and, most probably, there is no way
to give one definition of memory which would be good for all cases.
In different context the definition of memory can be different. For
example, considering limit theorems for partial sums the notions of
positive (long), zero (short), and negative memories are very
natural, while considering maximum operation instead of summation or
problem of large deviations, classification of memory properties can
be different - it seems that there is impossible to introduce
negative memory considering partial maxima operation, see
\cite{Samorodnitsky}, where the same value $1/\a$ serves as boundary
in the growth of partial maxima of stationary $S\a S$ random
variables.

The rest of the paper is organized as follows. In the second section
we consider $\a$-covariance for linear processes and fields with
infinite variance. The third section is devoted to $\a$-covariance
for random vectors defined as stochastic integrals and for stable
processes. In the fourth section we consider memory properties of
stationary random processes and fields. The last, fifth, section is,
may be, the most important, since we show that the notion  of
$\a$-covariance, which was introduced for $S\a S$ random vectors can
be  extended and generalized for much more broad class of random
vectors.

\section{Linear processes and fields}
\subsection{Linear processes}
We start with the case of linear processes. Let $\e_i, i\in \bz $
i.i.d. standard $S\a S$ random variables with ch.f. $\exp (-|t|^\a),
\ 0<\a \le 2$ (in the case of Gaussian distribution variance will be
not 1, but 2). We consider linear random process
\begin{equation}\label{linpr}
X_k=\sum_{j=0}^\infty c_j\e_{k-j}, \ k\in \bz,
\end{equation}
where $c_j, j\ge 0$ are real numbers satisfying condition
\begin{equation}\label{cond1}
A:=\sum_{j=0}^\infty |c_j|^\a <\infty.
\end{equation}
We get a stationary sequence of $S\a S$ random variables $X_k, \
k\in \bz$ with ch.f. $\exp (-A|t|^\a)$, and the main question is how
to measure the dependence between $X_0$ and $X_n$. Since bivariate
random vector $(X_0, X_n)$ is jointly $S\a S$, we can apply as a
measure of dependence $\a$-cc and $\a$-covariance. Let us denote
${\tilde \rho}_n={\tilde \rho} (X_0, X_n)$ and ${\rho}_n={\rho}
(X_0, X_n)$. To formulate our result we need some more notations.
Let
$$
a_{n,j}=(c_j, c_{n+j}), \ ||a_{n,j}||^2=(c_j^2+ c_{n+j}^2), \
{\tilde a}_{n,j}=(c_j, c_{n+j})||a_{n,j}||^{-1}
$$
$$
A_{1,n}=\sum_{j=0}^\infty \frac{c_j^2}{||a_{n,j}||^{2-\a}}, \
A_{2,n}=\sum_{j=0}^\infty \frac{c_{j+n}^2}{||a_{n,j}||^{2-\a}}, \
A_{n}=\sum_{j=0}^n |c_j|^\a.
$$
The convergence of the two above written series easily follows from
(\ref{cond1}), for example,
$$
A_{1,n}=\sum_{j=0}^\infty \frac{c_j^2}{(c_j^2+
c_{n+j}^2)^{(2-\a)/2}}\le \sum_{j=0}^\infty c_j^2|c_j|^{\a-2} =A.
$$
\begin{teo}\label{thm1} For a linear process $X_k$  from
(\ref{linpr}), satisfying (\ref{cond1}), we have
\begin{equation}\label{gapn}
{\tilde \rho}_n= \frac{\sum_{j=0}^\infty
c_jc_{j+n}||a_{n,j}||^{\a-2}}{\sqrt {A_{1,n}(A_{2,n}+A_{n-1})}},
\end{equation}
and
\begin{equation}\label{gapn1}
{ \rho}_n= \sum_{j=0}^\infty c_jc_{j+n}||a_{n,j}||^{\a-2}.
\end{equation}
\end{teo}

{\it Proof of Theorem \ref{thm1}. } We must find the spectral
measure of the $S\a S$ random vector $(X_0, X_n)= (\sum_{k=0}^\infty
c_k \e_{-k}, \ \sum_{k=0}^\infty c_k \e_{n-k})$. Denoting $\bt
=(t_1, t_2)$ and using the notations introduced before the
formulation of Theorem \ref{thm1}, we can write
\begin{eqnarray*}
 E\exp \{i(t_1X_0+t_2X_n)\} & = & E\exp \left
\{it_1\sum_{k=0}^\infty c_k \e_{-k} +it_2 \sum_{k=0}^\infty c_k
\e_{n-k})\right \}\nonumber \\
       & = & E\exp \left
\{\sum_{k=0}^\infty i (t_1c_k +t_2 c_{k+n}) \e_{-k} +
\sum_{k=0}^{n-1}it_2 c_k
\e_{n-k})\right \}\nonumber \\
 & = & \exp \left
\{-\left (\sum_{k=0}^\infty |(\bt ,{\tilde a}_{k,n}|^\a
||a_{n,k}||^\a+ \sum_{k=0}^{n-1}|c_k|^\a
|t_2|^\a\right ) \right \}\nonumber \\
\end{eqnarray*}
From this expression we see that the bivariate $S\a S$ random vector
$(X_0, X_n)$ has the symmetric spectral measure $\G_n$ concentrated
at points $(0, \pm 1), \pm {\tilde a}_{k,n}, \ k\ge 0$, namely,

\begin{equation}\label{gapn2}
\G_n (0, \pm 1)=\frac{1}{2}A_{n-1}, \quad \G_n (\pm {\tilde
a}_{k,n})=\frac{1}{2}||a_{n,k}||^\a.
\end{equation}
Due to (\ref{cond1}) this measure is finite:
$$
\G(S_2)=A_{n-1}+\sum_{k=0}^\infty ||a_{n,k}||^\a \le 2A.
$$
Having (\ref{gapn2}) we obtain (\ref{gapn}) and (\ref{gapn1}) by
simple calculations using definitions (\ref{gapdef}) and
(\ref{acovdef}). The theorem is proved. \halmos

 For two sequences $\{a_n\}$ and $\{b_n\}$ we shall write $a_n \sim
 b_n$,  if $\lim a_n b_n^{-1}=1$, and $a_n \simeq
 b_n$, if there exist two constants $0<K_1 <K_2<\infty$ such, that for
 sufficiently large $n$, \ $K_1\le a_n b_n^{-1}\le K_2$. We can
 formulate two simple properties of the introduced measures of
 dependence of a linear process.
 %Let ${\tilde X}_k$ be a linear
% process, defined as in (\ref{linpr}) only with a filter ${\tilde
% c}_j$ instead of $c_j.$

\begin{prop}\label{prop1} For any $c\neq 0$ we have
\begin{equation}\label{gapn3}
{\tilde \rho} (cX_0, cX_n)= {\tilde \rho} (X_0, X_n), \quad \rho
(cX_0, cX_n)=|c|^\a{\rho} (X_0, X_n).
\end{equation}
If ${\tilde c}_j \sim c_j$, then
\begin{equation}\label{gapn4}
\rho (X_0, X_n) \simeq \rho ({\tilde X}_0, {\tilde X}_n),
\end{equation}
where ${\tilde X}_n$ is defined by (\ref{linpr}), only with
coefficients $\{{\tilde c}_j\}$ instead of $\{c_j \}$.
 If, additionally to (\ref{cond1}), the following mild condition
\begin{equation}\label{cond3}
|c_{j+n}|\le k |c_j|, {\rm for \ \ all} \ j, n, \ {\rm and \ \ for \
\ some} \ k>0,
\end{equation}
is satisfied, then
\begin{equation}\label{equiv}
{\tilde \rho}_n \simeq {\rho}_n.
\end{equation}
\end{prop}
Due to the last property, as in the case of finite variance, we
shall deal mainly with $\a$-covariance, although there is a small
difference between these two cases with finite and infinite
variances: for a stationary sequence with finite variance,
correlation and covariance for all lags differs by a constant (equal
to the variance of the marginal distribution), while in the case of
a stationary sequence (\ref{linpr}) we have only (\ref{equiv}).

{\it Proof of Proposition \ref{prop1}.} The equalities (\ref{gapn3})
are obvious, since it is easy to see that if $\G_n$ and $\G_{n, c}$
are the spectral measures of  random vectors $(X_0, X_n)$ and
$(cX_0, cX_n),$ respectively, then $\G_{n, c}(ds)=|c|^\a
\G_{n}(ds).$

Let us denote ${\hat a}_{n,j}=({\tilde c}_j, {\tilde c}_{n+j})$.
Having relation ${\tilde c}_j= c_j (1+\d (j))$ with $\d(j)\to 0$,
for $j\to \infty $ and $|\d(j)|\le a$ for all $j$ and for a
sufficiently small $0<a<1$, we can get
$$
||{\hat a}_{n,j}||^2=|| a_{n,j}||^2(1+\d_1(j,n)),
$$
where $\d_1(j,n)\to 0$, for $j\to \infty $, uniformly with respect
to $n$, and $|\d_1(j, n)|\le a_1$. Here $a_1$ can be expressed by
$a$ and will be small for small $a$. Now we can write
\begin{equation}\label{gapn5}
\rho ({\tilde X}_0, {\tilde X}_n)=\sum_{j=0}^\infty {\tilde
c}_j{\tilde c}_{j+n}||{\hat a}_{n,j}||^{\a-2}=\sum_{j=0}^\infty
c_jc_{j+n}||a_{n,j}||^{\a-2} (1+\d_2(j,n)),
\end{equation}
where $\d_2$ is obtained from the formal relation
$$
1+\d_2(j,n)=\frac{(1+\d (j))(1+\d (j+n))}{(1+\d_1(j,n))^{(2-\a)/2}}.
$$
Again, it can be shown that $\d_2$ has the same properties as $\d_1$
and can be bounded $|\d_2(j, n)|\le a_2<1$ if $a$ is chosen
sufficiently small. The relation (\ref{gapn4}) follows from
(\ref{gapn5}).

To prove (\ref{equiv}) we need to show that the quantity
$A_{1,n}(A_{2,n}+A_{n-1})$ is bounded from above and below by some
constants. The bound from above
$$
A_{1,n}(A_{2,n}+A_{n-1})\le A^2
$$
is easy, while for the lower bound for $A_{1,n}$ we use
(\ref{cond3}):
$$
A_{1,n}\ge \sum_{j=0}^\infty
\frac{c_j^2}{(c_j^2(1+k^2))^{(2-\a)/2}}=(1+k^2))^{(\a-2)/2}A.
$$
For any $\epsilon >0$,  we can find $N$ such that $\sum_{j=n}^\infty
|c_j|^\a <\epsilon$ for all $n>N,$, therefore
$$
A_{2,n}+A_{n-1}\ge A_{n-1}=A-\sum_{j=n}^\infty |c_j|^\a \ge
(1-\epsilon)A.
$$
From these estimates the relation (\ref{equiv}) follows, and the
proposition is proved.  \halmos
 In special cases of the filter of a linear process
(\ref{linpr}) we have the following corollaries. %In the case of
%exponentially decaying coefficients of the filter we have very
%simple relation.
\begin{cor}\label{corol1} Let $c_j\sim 2^{-j}$ then $ \rho_n\sim
C(\a)2^{-n}$.

If $c_j\sim j^{-\b}, \ \b>1/\a$, then: in the case $0<\a \le 1$
$$
 \rho_n\simeq C(\a, \b)n^{1-\b \a},
$$
in the case $1<\a \le 2$
$$
 \rho_n\simeq  \left \{ \begin{array}{ll}
                           C(\a, \b)n^{1-\b\a}, & \mbox{if  $1/\a <\b
                           <1/(\a-1)$,}\\
                           C(\a)n^{-\b}(1+\bbI (\b=1/(\a-1))\ln n), & \mbox{if
                           $\b \ge
                           1/(\a-1)$},
                           \end{array}
                      \right.
$$
where, as usual, $\bbI (A)$ stands for the indicator function of an
event $A$.
\end{cor}
{\it Proof of Corollary \ref{corol1}.} Taking $c_j=2^{-j}$ we simply
have
$$
\rho_n= \sum_{j=0}^\infty
2^{-j}2^{-j-n}(2^{-j}(1+2^{-2n})^{1/2})^{\a-2}=2^{-n}\frac{(1+2^{-2n})^{(\a-2)/2}}
{1-2^{-\a }}.
$$
 Now let us  take $c_j=j^{-\b}, \ j\ge
1, \ c_0=1$ and in the definition (\ref{gapn1}) of $\rho_n$ we
separate two first terms
\begin{equation}\label{gapn6}
c_0c_{n}||a_{n,0}||^{\a-2}+ c_1c_{n+1}||a_{n,1}||^{\a-2}\sim
n^{-\b},
\end{equation}
then we can write
$$
I_n :=\sum_{j=2}^\infty c_jc_{j+n}||a_{n,j}||^{\a-2}\simeq
\int_1^\infty \frac{x^{-\b}(x+n)^{-\b} dx}{\left
(x^{-2\b}+(x+n)^{-2\b} \right )^{(2-\a)/2}}.
$$
After  change of variables we get
\begin{equation}\label{gapn7}
I_n=n^{1-\b \a}(I_{n,1}+I_{n,2}),
\end{equation}
where
$$
I_{n,1}=\int_{1/n}^1 \frac{y^{-\b}(1+y)^{-\b} dy}{\left
(y^{-2\b}+(1+y)^{-2\b} \right )^{(2-\a)/2}},
$$
and $I_{n,2}$ is the integral of the same function over interval
$(1, \infty).$ Therefore, $I_{n,2}$ is independent of $n$, and due
to the condition $\b \a>1$, is a constant depending on $\a$ and $\b$
only:
$$
I_{n,2}=\int_1^\infty \frac{dy}{y^{\b \a}g(y)} =C(\a, \b).
$$
Here $g(y)$ is a function, bounded by positive constants from below
and above for all $1\le y <\infty$. Now we estimate $I_{n,1}$.
Again, it is easy to see that
$$
I_{n,1}=\int_{1/n}^1 \frac{dy}{y^{\b(\a-1)}h(y)},
$$
where $1\le h(y)\le 2^{\b}(1+2^{-2\b})^{(2-\a)/2}$ for $0\le y\le
1$. If $\a \le 1$ then $\a-1\le 0$, and we get that $I_{n,1}$ is a
constant, depending on $\a$ and $\b$. In the case of $1<\a \le 2$
and $1/\a <\b <1/(\a-1)$, we have again that $I_{n,1}$ is a
constant, while if $\b=1/(\a-1)$ we get that $I_{n,1}$ is of the
order $\ln n$. Finally, if $\b>1/(\a-1)$, then $I_{n,1}\simeq C(\a,
\b)n^{\b(\a-1)-1}$, and we get from (\ref{gapn7}) that $I_n$ is of
the order $n^{-\b}$, as the first two terms in (\ref{gapn6}). The
corollary is proved. \halmos

Now we can compare our result for $\a$-covariance of a linear
process (\ref{linpr}) with other measures of dependence.  There were
several papers dealing with measures of dependence of linear
processes with innovations with infinite variance, see, for example
\cite{Kokoszka1}, \cite{Kokoszka}, \cite{Kokoszka2} and references
there. In these papers the expressions of the covariation and the
codifference for the process (\ref{linpr}) were given:
$$
\t (X_0, X_n)=\sum_{j=0}^\infty (|c_j|^\a +|c_{j+n}|^\a
-|c_j-c_{j+n}|^\a),
$$
$$
[X_0, X_n]_\a =\sum_{j=0}^\infty  c_{j+n}c_j^{\langle \a-1\rangle}.
$$
The asymptotic of these quantities was investigated for
$FARIMA(p,q,d)$ process  (see (2.2), (2.6) and (2.7) formulas in
\cite{Kokoszka1}), which is of the form (\ref{linpr}) with specific
coefficients $c_j, \ j\ge 0$. Namely, in \cite{Kokoszka1} it was
shown that these coefficients satisfy the following relation (see
Lemma 3.2 and Corollary 3.1 there): if $\a(d-1)<-1$, then
\begin{equation}\label{farimacoef}
c_j =C(p,q,d)j^{d-1}(1+O(j^{-1})).
\end{equation}
Quantity $d$ here and $\b$ used in Corollary \ref{corol1} are
related by equality $\b=1-d$. In \cite{Kokoszka2} general case of
(\ref{linpr}) is investigated under conditions which are slightly
different from those used in Corollary \ref{corol1}: the
coefficients of the filter satisfy more general condition $c_j=U(j)$
where $U$ is regularly varying function with the index $-\b$, but
there are some conditions of the type (\ref{farimacoef}). Therefore
in \cite{Kokoszka2} the  asymptotic relation $\sim$ is obtained,
while we have weaker relation $\simeq$, but the order of decay of
the quantity $\t_n:=\t (X_0, X_n)$ is the same as of the
$\a$-covariance in Corollary \ref{corol1}. We can mention that in
both papers  \cite{Kokoszka1} and \cite{Kokoszka2} the case
$\b=1/(\a-1)$ is excluded from formulation, only mentioning that
there is "phase transition" (see the remark before Theorem 4.1 in
\cite{Kokoszka1}). Also it is interesting to note that in the case
of exponentially decreasing filter ($c_i\sim 2^{-j}$) the
codifference is decreasing as $2^{-\a n}$ for $0<\a <1$, while from
Corollary \ref{corol1}we have $ \rho_n\simeq C(\a)2^{-n}$ for all
$0<\a \le 2$

\subsection{Linear fields}
%\bigskip

As it was mentioned in the introduction,  $\a$-covariance, as the
measure of dependence, can be easily applied for linear random
fields  on $\bz^d$ with $d\ge 2$ and $S\a S$ innovations. But since
the notation and formulations become more complicated we restrict
ourselves to the case $d=2$ and formulation of expression of
$\a$-covariance via coefficients of a filter of a random field.
 Let $\e_{i,j}, \ (i, j) \in \bz^2 $ be
i.i.d. standard $S\a S$ random variables with ch.f. $\exp (-|t|^\a),
\ 0<\a \le 2$. Now we consider linear random field
\begin{equation}\label{linfield}
X_{k,l}=\sum_{i,j=0}^\infty c_{i,j}\e_{k-i, l-j}, \ (k, l)\in \bz^2,
\end{equation}
where $c_{i, j}, \ i\ge 1,  j\ge 0,$ are real numbers satisfying
condition
\begin{equation}\label{cond2}
A_1:=\sum_{i, j=0}^\infty |c_{i,j}|^\a <\infty.
\end{equation}
We are interested how strongly dependent are random variables
$X_{0,0}$ and $X_{n,m}$. Let us denote ${\tilde \rho}_{n,m}={\tilde
\rho} (X_{0,0}, X_{n,m}), \quad \rho_{n,m}=\rho (X_{0,0}, X_{n,m})$.
We shall consider two cases: $n>0, m>0$ and $n>0, m<0$, since due to
the stationarity the remaining two cases can be transformed into the
previous, for example, $(X_{0,0}, X_{n,m})$  has the same
distribution as $(X_{-n,-m}, X_{0,0})$ and
$\rho_{n,m}=\rho_{|n|,|m|}$ for $n<0, m<0$. Similarly,
$\rho_{n,m}=\rho_{-n,-m},$ for $n>0, m<0$.   Comparing with the case
of linear processes now we need more complicated notations. First,
let us consider the case $n>0, m>0$. Let us denote
$$
a_{i,j}^{(n,m)}=(c_{i,j}, c_{i+n,j+m}), \
||a_{i,j}^{(n,m)}||^2=c_{i,j}^2+ c_{i+n,j+m}^2,
$$
$$ {\tilde a}_{i,j}^{(n,m)}=(c_{i,j},
c_{i+n,j+m})||a_{i,j}^{(n,m)}||^{-1}
$$
Denote the following four regions of $\bz_{+}^2$:
$$
I_1=\{(i,j): 0\le i \le n-1, \ 0\le j \le m-1, \}, \ I_2=\{(i,j):
0\le i \le n-1, \  j \ge m, \},
$$
$$ I_3=\{(i,j): i \ge n, \ 0\le j \le m-1, \}, \ I_4=\{(i,j):  i \ge n, \ j \ge m,
\}.
$$
Also denote $\sum_k= \sum \sum_{(i,j)\in I_k}, \ k=1,2,3,4$, and
$\sum_0= \sum \sum_{(i,j)\in \bz_{+}^2}$. Then we define
$$
A_{1, n, m}:=\left ( {\sum}_1 +{\sum}_2 +{\sum}_3\right
)|c_{i,j}|^\a, \quad A_{2, n, m}:={\sum}_0
c_{i,j}^2||a_{i,j}^{(n,m)}||^{\a-2}
$$
$$
A_{3, n, m}:={\sum}_0 c_{i+n,j+m}^2||a_{i,j}^{(n,m)}||^{\a-2}
$$

In the case $n>0, m<0$ we need the following notation:
$$
J_1=\{(i,j): 0\le i <\infty, \ 0\le j \le |m|-1, \}, \ J_2=\{(i,j):
0\le i <\infty, \  j \ge 0, \},
$$
%$$
%J_3=\{(i,j): -n\le i \le -1, \ |m|\le j <\infty, \},
%$$
and $\sum^{(k)}= \sum \sum_{(i,j)\in J_k}, \ k=1,2.$ Then we define
$$
B_{n,m}^{(i)}={\sum}^{(i)}|c_{i,j}|^\a, \quad i=1,2,
$$
$$
B_{n,m}^{(3)}={\sum}_0 c_{i,j+|m|}^2(c_{i,j+|m|}^2
+c_{i+n,j}^2)^{(\a-2)/2}, \quad B_{n,m}^{(4)}={\sum}_0
c_{i,j+|m|}^2(c_{i,j+|m|}^2 +c_{i+n,j}^2)^{(\a-2)/2}.
$$

Now we are able to formulate our result for linear random field
(\ref{linfield}).
\begin{teo}\label{thm2} For a linear field $X_{k,l}$  from
(\ref{linfield}), satisfying (\ref{cond2}), for $n>0, m>0$, we have
$$
{\tilde \rho}_{n,m}=\frac{{\sum}_0 c_{i,j}c_{i+n,j+m}
||a_{i,j}^{(n,m)}||^{\a-2}}{\sqrt {A_{2,n,m}(A_{1,n,m}+A_{3,n,m})}},
$$
and
$$
 \rho_{n,m} = {\sum}_0 c_{i,j}c_{i+n,j+m}
||a_{i,j}^{(n,m)}||^{\a-2}.% ,\quad {\rm as} \ \min (n,m)\to \infty.
$$
 If $n>0, m<0,$  then
$$
{\tilde \rho}_{n,m}=\frac{{\sum}_0 c_{i,j+|m|}c_{i+n,j}
(c_{i,j+|m|}^2 +c_{i+n,j}^2)^{(\a-2)/2}}{\sqrt
{(B_{n,m}^{(1)}+B_{n,m}^{(3)})(B_{n,m}^{(2)}+B_{n,m}^{(4)})}},
$$
and
$$
 \rho_{n,m} ={\sum}_0 c_{i,j+|m|}c_{i+n,j}
(c_{i,j+|m|}^2 +c_{i+n,j}^2)^{(\a-2)/2}.
$$

\end{teo}

Clearly, for linear fields we can formulate the same properties of
${\tilde \rho}_{n,m}$ and $ \rho_{n,m}$ as in Proposition
\ref{prop1} for linear processes. Also we can easily calculate these
quantities for filters with regular behavior. As an example we
provide one such result. Suppose that coefficients of the filter are
hyperbolically decaying: $c_{i, j}\sim i^{-\b_1}j^{-\b_2}$ with
 $\b_k>1/\a, \ k=1,2$. For the convenience let us denote ${\bar \b}=(\b_1, \b_2).$
\begin{cor}\label{corol4} Let $c_{i, j}\sim i^{-\b_1}j^{-\b_2}$ with
 $\b_k>1/\a, \ k=1,2$,  and $m>0, \ n>0.$ Then, if $0<\a \le 1,$
$$
 \rho_{n, m}\simeq C(\a,{\bar \b})n^{1-\b_1 \a}m^{1-\b_2 \a}, \ \b_i>1/\a;
$$
if $1<\a \le 2,$ then
$$
 \rho_{n, m}\simeq  \left \{ \begin{array}{ll}
                           C(\a, {\bar \b})n^{1-\b_1 \a}m^{1-\b_2 \a}, & \mbox{if  $1/\a
                           <\b_i<(\a-1)^{-1}, i=1,2 $,}\\
                           C(\a, {\bar \b})n^{-\b_1}m^{-\b_2}, & \mbox{if  $\b_i>
                           1/(\a-1), i=1,2$.}
                           \end{array}
                      \right.
$$
If $\b_1\ge 1/(\a-1)$,  $1/\a <\b_2<(\a-1)^{-1} $, then
$$
\rho_{n, m}\simeq C(\a, {\bar \b})n^{-\b_1}m^{1-\b_2 \a}(1+\bbI
(\b_1=(\a-1)^{-1})\ln n),
$$
and if $\b_2\ge 1/(\a-1)$, $1/\a <\b_1<(\a-1)^{-1} $, then
$$
\rho_{n, m}\simeq C(\a, {\bar \b})n^{1-\b_1\a}m^{-\b_2}(1+\bbI
(\b_2=(\a-1)^{-1})\ln n).
$$
\end{cor}
This result is in accordance with the result of Corollary
\ref{corol1}. We omit the calculations needed to prove these
relations, since they are very similar to those used for linear
processes.

\vspace{10pt}
\section{$\a$-covariance for stochastic integrals and stable processes}

It is well-known what important role in the theory of stable random
vectors and processes play $\a$-stable stochastic integrals, that
is, integrals of non-random functions with respect to an $\a$-stable
random measures. The biggest part of the monograph \cite{Samorod} is
devoted to these integrals, therefore we do not provide definitions
of these notions (but we shall try to keep the same notation as in
\cite{Samorod}), refereing a reader to this monograph. Let $(E,
{\cal E}, m)$ be a measurable space with a $\s$-finite measure $m$,
and let $M$ be an $S\a S$ random measure, that is, we take  the
so-called skewness intensity function $\b (x)\equiv 0$ in general
Definition 3.3.1 in \cite{Samorod}  of $\a $-stable random measure.
This is done for the reason that $\a$-covariance we defined (till
now, see Section 5 for extension of definition of $\a$-covariance)
only for $S\a S$ random vectors. Taking $f\in L^{\a}(E, {\cal E},
m)$, we get a $S\a S$ random variable
$$
X=\int_E f(x)M(dx),
$$
while taking a collection $f_i\in L^{\a}(E, {\cal E}, m), \ i=1,
\dots ,k$, \ we get a $S\a S$ random vector
$$
(X_1, \dots ,X_k), \quad X_i=\int_E f_i(x)M(dx).
$$
Taking a family of functions  $\{f_t, t \in T\} \subset L^{\a}(E,
{\cal E}, m)$ we get a $S\a S$ random process
$$
X(t)=\int_E f_t(x)M(dx), \ t\in T.
$$
Many well-known $S\a S$ random processes are obtained in this way.
Namely, a moving average process is obtained with $E=\R$,
$m$=Lebesgue measure, and $f_t(x)=f(t-x)$:
\begin{equation}\label{movaver}
X(t)=\int_{\R} f(t-x)M(dx), \ t\in \R.
\end{equation}
An Ornstein-Uhlenbeck process
\begin{equation}\label{ornstein}
X(t)=\int_{-\infty}^t \exp \{-\l(t-x)\}M(dx), \ t\in \R,
\end{equation}
is obtained from (\ref{movaver}) by taking $f(x)=\exp \{-\l x\}\bbI
(x\ge 0)$. In a similar way, i.e., by choosing appropriate family of
functions $f_t(x)$, we can get a symmetric linear fractional stable
motion, a log-fractional stable motion (see Ch. 3.6 in
\cite{Samorod}). It is worth to mention that linear processes and
fields, considered above can be obtained in the same way from
(\ref{movaver}), taking $E=\bz$ (or $E=\bz^2$ in the case of
fields),  $m$ as a counting measure, and $f(k)=c_k \bbI (k\ge 0)$.
The covariation and the codifference were defined for these more
general objects, and they were extensively studied during last three
decades. Most of these results are given in \cite{Samorod}, see Ch.
4.7 therein, where the codifference function is calculated for many
stationary $S\a S$ processes. Our goal is to show that
$\a$-covariance is equally good measure of dependence, more easily
dealt with,  and  in some cases even better reflect dependence.

Let $(X_1, X_2)$ be a bivariate $S\a S$ random vector, defined by
means of stochastic integrals, i.e.,
$$
(X_1, X_2)\lygpas \left (\int_E f_1(x)M(dx), \ \int_E f_2(x)M(dx),
\right )
$$
where $\lygpas $ stands for equality in distribution. In Ch. 3.2 in
\cite{Samorod} it is shown how to express spectral measure $\G$ of
the random vector $(X_1, X_2)$ via control measure $m$ and functions
$f_i$,  also  the expressions of the covariation (in the case $1<\a
\le 2$) and the codifference are given:
$$
[X_1, X_2]_\a=\int_E f_1(x)f_2(x)^{\langle \a-1\rangle}m(dx),
$$
$$
\t (X_1, X_2)=\int_E \left (|f_1(x)|^\a
+|f_2(x)|^\a-|f_1(x)-f_2(x)|^\a\right )m(dx).
$$
It is not difficult to write the expressions of $\a$-covariance and
$\a$-correlation in these terms:
\begin{equation}\label{acovardef1}
\rho (X_1, X_2)=\int_E f_1(x)f_2(x)||{\bar f}(x)||^{\a -2}m(dx),
\end{equation}

$$
{\tilde \rho}(X_1, X_2)=\frac{\int_E f_1(x)f_2(x)||{\bar f}(x)||^{\a
-2}m(dx)}{\left (\int_E\frac{ f_1^2(x)}{||{\bar f}(x)||^{2-\a
}}m(dx)\int_E \frac{f_2^2(x)}{||{\bar f}(x)||^{2-\a }}m(dx) \right
)^{1/2}},
$$
where $||{\bar f}(x)||=\left (f_1^2(x) +f_2^2(x) \right )^{1/2}.$
Formally in the above written formulae we should integrate over
$E_+= \{ x\in L^\a: ||{\bar f}(x)||>0 \}$, but here we agree that
integrand is equal to zero if $||{\bar f}(x)||=0.$ Comparing
expressions of the codifference and $\a$-covariance, we see that the
integrand in (\ref{acovardef1}) is more simple to deal with. We
shall demonstrate this by calculating $\a$-covariance function for
Ornstein-Uhlenbeck process (\ref{ornstein}). Let us denote $\rho
(t)=\rho (X(0), X(t))$ and the normalized $\a$-covariance function
${\bar \rho}(t):=\rho (t)(\rho (0))^{-1} $. Since $\rho (-t)=\rho
(t)$, for $t>0$, it is sufficient to consider the case $t>0.$
\begin{prop}\label{ornsteinprop}
Let $X$ be the process defined in (\ref{ornstein}). For $t>0$ we
have
\begin{equation}\label{ornstein1}
\rho (t)=\frac {1}{\a \l (1+\exp (-2\l t))^{(2-\a)/2}}e^{-\l t}
\end{equation}
and
\begin{equation}\label{ornstein2}
{\bar \rho} (t)=\frac {1}{ (2^{-1}(1+\exp (-2\l
t)))^{(2-\a)/2}}e^{-\l t}.
\end{equation}
Thus, as $t\to \infty$,
\begin{equation}\label{ornstein3}
{\bar \rho} (t)  \sim 2^{(2-\a)/2}e^{-\l t}.
\end{equation}
\end{prop}
{\it Proof of Proposition \ref{ornsteinprop}.}  From
(\ref{ornstein}) we see that we must apply (\ref{acovardef1}) with
$$
f_1(x)=\exp (\l x)\bbI (x\le 0), \quad f_2(x)=\exp (\l x -\l t)\bbI
(x\le t).
$$
Then
$$
\rho (t)= \int_{-\infty}^\infty \frac{\exp (\l x)\bbI (x\le 0)\exp
(\l x -\l t)\bbI (x\le t)}{\left (\exp (2\l x)\bbI (x\le 0)+\exp
(2\l x -2\l t)\bbI (x\le t)\right )^{(2-\a)/2}} dx,
$$
and simple integration gives us (\ref{ornstein1}). Since $\rho
(0)=(\a \l)^{-1}2^{(\a-2)/2},$ the equality (\ref{ornstein2}) is
obtained from (\ref{ornstein1}), and the relation (\ref{ornstein3})
is obvious. \halmos

We can compare these results with corresponding results for the
codifference, provided in Example 4.7.1 in \cite{Samorod}. If we
denote by $\t(t)=\t(X(0), X(t))$ and normalized the codifference
function  by ${\bar \t}(t)=\t (t)(||X(0)||_\a^\a)^{-1}$, where
$||X(0)||_\a^\a$ stands for the scale parameter of $S\a S$ random
variable $X(0),$ then
\begin{equation}\label{ornsteincodiff}
\t (t)=\frac{1}{\a \l}\left (1-(1-\exp (-\l t))^{\a}+e^{-\a \l
t}\right )
\end{equation}
and
\begin{equation}\label{ornsteincodiff1}
 {\bar \t}(t) \sim  \left \{ \begin{array}{ll}
                           \a \exp (-\l t), & \mbox{if  $1<\a <2$,}\\
                           2 \exp (-\l t), & \mbox{if $\a=1$}, \\
                           \exp (-\a \l t), & \mbox{if $0<\a <1$}.
                           \end{array}
                      \right.
\end{equation}
Comparing (\ref{ornstein1}) with (\ref{ornsteincodiff}) we see that
expression  for $\a$-covariance is more simple and, the most
important, gives exponential decay independent of $\a$, while for
the codifference in (\ref{ornsteincodiff1}), in the range $0<\a<1$,
there is $\a$ in the exponent. We see the same effect as in the case
of linear processes with exponentially decaying filters, see the
discussion at the end of subsection 2.1. The constant in the
asymptotic of the normalized $\a$-covariance function in
(\ref{ornstein3}) is continuous in $\a$ and varies in the small
interval $(1/2, 1]$, while in (\ref{ornsteincodiff1}) dependence of
the constant on $\a$ is discontinuous  at $\a=1$. Dependence of the
$\a$-covariance function on the parameter $\l$ is the same as of the
codifference function, namely, if we denote by $\t (t, \l)$ and
$\rho (t, \l)$ the codifference  and $\a$-covariance functions,
respectively,  of the Ornstein-Uhlenbeck process $X(t)$ with
parameter $\l$, then we have
\begin{equation}\label{ornsteincodiff2}
\t(t, \l_1)< \t(t, \l_2), \quad \rho (t, \l_1)< \rho (t, \l_2),
\end{equation}
for $\l_2 <\l_1$ and for all $0<\a\le 2.$ To prove the second
relation in (\ref{ornsteincodiff2}) (the first one is proved in
\cite{Samorod}, see p. 210 therein) we consider (for a fixed $t>0$)
the function
$$
h(\l)=\frac{\exp (-\l t)}{\l (1+\exp (-2\l t))^{(2-\a)/2}},
$$
and it is easy to show that $h'(\l)<0$ for all $\l >0, \ t>0, \ 0,\a
\le 2$.
 Here it is worth to mention that in the case $\a=2$
the constants obtained from formulas (\ref{ornstein1}) and
(\ref{ornstein3}) do not coincide with constants, given in
\cite{Samorod}, see p. 210:  if $\a=2$, then  $\rho (t)=\t (t)/2$
(this equality can be seen also from the relation $|s_1|^2 +|s_2|^2
-|s_1-s_2|^2=2s_1 s_2$), and this difference comes from the fact
that characteristic function of the standard $S\a S$ random variable
is $\exp (-|t|^\a)$, while for Gaussian standard random variable
this function is $\exp (-t^2/2)$.
\medskip

We took the Ornstein-Uhlenbeck process as an example from large
class of processes, whose finite dimensional distributions are $S\a
S$. In the theory of stochastic processes there are important
classes of $\a$-stable processes, such as sub-Gaussian, moving
averages, harmonizable processes, fractional stable noises, etc.,
rather detailed study of such processes is presented in the
fundamental monograph \cite{Samorod}. As the main tool in
\cite{Samorod} to study dependence for these processes is used the
codifference function. We believe (and this belief is based on the
extensions and generalizations given in the last section) that in
the case of infinite variance $\a$-covariance function is a better
substitute for usual covariance function, although we admit that  a
lot of work must be done - during thirty years there were a lot of
papers devoted to the codifference and covariation, while this paper
is the first one after 1976 paper \cite{Paul3} (where this notion
was only introduced) devoted to $\a$-covariance function.

\section{Short, long and negative memories}

\subsection{Memory properties for random  processes}

The importance of  notions of long-range and short-range dependence
and notion of  memory   in the theory of stochastic processes and
fields and, in particular, in time series analysis is well-known.
The number of monographs and papers devoted to these notions are
growing steadily, and this can be explained from one hand, by usage
of these notions in many areas, ranging from econometrics and
finance to hydrology and climate studies, on the other hand, by the
complexity of these notions , complicated relations with other
notions. We refer to important survey paper \cite{Samorodnitsky1}
and recent monograph \cite{Giraitis} which gave impetus to look at
these notions for sequences with infinite variance.

If we consider a stationary sequence with finite variance, there are
several ways to define long-range dependence, in \cite{Guegan} there
are provided even 8 definitions, but the  three of them are main
(other 5 are only modifications): via covariance function; via
spectral density; and via the growth of partial sums (the so-called
Allen variance). Let us note that in many papers the notions
"long-range dependence" and "long memory" are used as synonyms. In
the above cited  monograph \cite{Giraitis} in the subject list there
is no notions "long-range dependence" or "short-range dependence"
and stationary processes  with finite second moment are divided into
processes with short, long and negative memory; this is done by
means of the covariance function. Namely, the following definition
is given in \cite{Giraitis}, see Definition $3.1.2$ there.

\begin{definition}\label{def1} A covariance stationary mean zero process $Y_t, \ t\in \bz,$
with covariance function $\g_k=EY_0Y_k$ has: long memory if
$\sum_{k\in \bz} |\g_n|=\infty;$ short memory if $\sum_{k\in \bz}
|\g_n|<\infty \quad {\rm and} \ \ \sum_{j\in \bz}^\infty \g_j > 0;$
and negative memory (or antipersistence) if $\sum_{k\in \bz}
|\g_n|<\infty \quad {\rm and} \ \ \sum_{j\in \bz}^\infty \g_j = 0.$
\end{definition}

One of the main application of this definition is to linear
processes with white-noise innovations (defined or by the relation
(\ref{linpr}), either by the analogous relation with summation over
all $\bz$), and this is due to the rather simple expression of
$\g_j$ via coefficients of the filter and the relation $\sum_{j\in
\bz}^\infty \g_j =(\sum_{j\in \bz} c_j)^2$. This relation allows to
use the sum $\sum_{j\in \bz} c_j$ for separation of short and
 negative memories. Therefore, it
 is tempting by means of the $\a$-covariance to define the
same notions for linear processes (\ref{linpr}), namely, we would
say that the process (\ref{linpr}) has: short memory if $ \sum_{k\in
\bz} |\rho_n|<\infty \quad {\rm and} \ \ \sum_{j=0}^\infty c_j \ne
0;$ long memory if $ \sum_{k\in \bz} |\rho_n|=\infty;$ and negative
memory if $ \sum_{k\in \bz} |\rho_n|<\infty \quad {\rm and} \ \
\sum_{j=0}^\infty c_j = 0.$ Unfortunately, such classification of
linear processes is unappropriate and useless for the following
reason. We know that dependence and, particularly, memory properties
play an important role in establishing limit properties of partial
sum processes constructed from stationary sequences under
consideration. These relations between memory properties and limit
theorems for partial sum processes are deeply investigated, and we
refer to the monograph \cite{Giraitis} where these relations are
given in details. Here we shall mention only that if $S_n$ stands
for the partial sum of a linear process with innovations with a
finite variance and regularly varying filter $c_k\sim k^{-1+d}$,
then $Var S_n \sim n$ in the case of short memory, and $Var S_n \sim
n^{1+2d}$ in the case of long memory $(0<d<1/2)$ and in the case of
negative memory ($-1/2<d<0$ and $\sum_{j=0}^\infty c_j = 0$). The
similar situation is in the case of linear processes with infinite
variance innovations. In \cite{Astrauskas} general limit theorems
for the partial sum process formed by a linear process with
innovations belonging to the domain of attraction of a stable law
are proved, and these limit theorems can be taken as a basis for
classification of linear processes with respect to memory
properties. We shall provide here simplified version of the results
from \cite{Astrauskas} avoiding more complicated formulations
involving slowly varying functions. Let us consider a linear process
\begin{equation}\label{linpr2}
X_k=\sum_{j=0}^\infty c_j\eta_{k-j}, \ k\in \bz,
\end{equation}
where $\{\eta_i, i\in \bz \},$ are i.i.d. random variables belonging
to the normal domain of attraction of a standard $S\a S$ random
variable with ch.f. $\exp (-|t|^\a), \ 0<\a < 2$, and a filter
$\{c_j, \ j\ge 0\} $ satisfies the relation $|c_j|\sim j^{-\b}$. Let
us consider the convergence of finite-dimensional distributions of
the process
\begin{equation}\label{normlinpr2}
Y_n(t)=A_n^{-1}\sum_{k=1}^{[nt]}X_k.
\end{equation}
In \cite{Astrauskas} three cases are separated.

(i) If  $\sum_{j} |c_j|<\infty$  and $\sum_{j} c_j \ne 0$, then
$A_n$ grows as $n^{1/\a}$, and the limit process is $\a$-stable L{\'
e}vy motion.

(ii) If $\a>1$ and $1/\a<\b<1$, then $A_n$ grows as $n^{1/\a+1-\b}$
(more rapidly comparing with the case (i)) and the limit process is
a linear fractional stable  motion.

(iii) Let $0<\a<2, \ \max (1, 1/\a)<\b<1+1/\a,$ and
\begin{equation}\label{cond4}
\sum_{j=0}^n c_j\sim (\b-1)^{-1}n^{1-\b},
\end{equation}
than $A_n$ grows as $n^{1/\a+1-\b}$ (now more slowly comparing with
the case (i)) and again the limit process is a linear fractional
stable motion.

It is necessary to note, that condition (\ref{cond4}) is stronger
then condition $\sum_{j=0}^\infty c_j=0$. Also from these results
and Corollary \ref{corol1}  it is clear that in the case $\a<2$
memory properties can not be characterized by the convergence or
divergence of series $\sum_n \rho_n$, as it was proposed above.
Therefore, it seems more natural in the case $\a<2$ memory
properties to define according the growth of normalizing constants
in limit theorems for partial sums, and this can be done not only
for linear processes but for general stationary sequences.

Let $\{\xi_i, \ i\in \bz\}$ be a strictly stationary sequence which
is jointly regularly varying with the index $0<\a<2$ (for the
definition of jointly regularly varying sequence see, for example
\cite{Basrak}). In order not to deal with centering we additionally
assume that $E\xi_0=0$ if $\a>1$ and that $\xi_0$ is symmetric if
$\a=1.$ Let us denote
$$
S_n(t)= \sum_{k=1}^{[nt]}\xi_k,
$$
and we suppose that there exists a sequence of normalizing constants
$A_n$ such that finite-dimensional distributions (f.d.d) of the
process $A_n^{-1}S_n(t)$ converges weakly to corresponding f.d.d. of
some stable processes (in particular, distribution of
$A_n^{-1}S_n(1)$ converges to an $\a$-stable law).

\begin{definition}\label{def2}We say that the sequence $\{\xi_i, \ i\in \bz\}$ is:  of short
memory if $A_n=n^{1/\a}L(n)$ with some slowly varying function $L$
and the limit process is the L{\' e}vy stable motion; of long memory
if $A_n=n^{1/\a+\d}L(n)$ with some $0<\d<1-1/\a$ and the limit
process is a linear fractional stable motion; of negative memory if
$A_n=n^{1/\a+\d}L(n)$ with some $-1/\a<\d<0$ and the limit process
is a fractional stable motion.
\end{definition}
It is worth to note that in this definition the condition
$0<\d<1-1/\a$ means that long memory can be only in the case $\a>1$.
Heuristically it can be explained as follows: in the case $0<\a<1$
and independent $\{\xi_i, \ i\in \bz\}$ the normalizing sequence
satisfies $A_n^\a \sim n$, and it is clear that for any stationary
sequence the scale parameter of $S_n$ can not grow faster due to the
moment inequality $E|\sum_{k=1}^{n}\xi_k|^\b\le nE|\xi_1|^\b$ for
any $\b<\a<1$.

 In the case of linear processes the above given three cases from \cite{Astrauskas}, formulated above,
exactly gives us three cases of memory, defined in Definition
\ref{def2}, namely,

 (i) we have long memory if $\a>1, \ 1/\a
 <\b<1$ ($A_n\sim n^{1/\a+\d}$ with $0<\d=1-\b<1-1/\a$),

 (ii) short memory if $0<\a<2, \ \max (1/\a, 1)
 <\b$ and $\sum_j c_j\ne 0$, ($A_n\sim n^{1/\a}$),

 (iii)  negative memory if $0<\a<2, \ \max (1/\a, 1)
 <\b<1+1/\a$ and $\sum_j c_j =0$ ($A_n\sim n^{1/\a+\d}$ with
 $-1/\a<\d=1-\b<0$).

\smallskip

Thus we have two different definitions of long, short and negative
memories in cases $0<\a<2$ and  $\a=2$ with finite variance. Of
course, the case $\a=2$ with infinite variance can be included
without any difficulties into Definition \ref{def2} with obvious
changes of limits in the definition in this case. We claim that more
logical definition is the second one, thus, the memory properties in
the case of stationary sequences with finite variance should be
defined as in Definition \ref{def2}. To justify this claim we shall
provide several simple examples. Let us  consider simple linear
process (\ref{linpr2}) with $E\eta_1^2<\infty$ and $|c_k|=k^{-\b}, \
k\ge 1$ ($c_0$ we shall define in several ways). For such simple
model covariances and normalization constants are easily calculated,
and our goal is to show that  Definition \ref{def2} is more logical.
 We consider the growth of the variance
$$
A_n^2= Var \sum_{k=1}^{n}X_k.
$$
In \cite{Giraitis} the asymptotic of this variance is obtained
investigating behavior of covariances $\g_k$, but for our purpose it
is more convenient to write explicit expression of $A_n^2$ via
coefficients $c_k$, namely
\begin{equation}\label{varsn}
A_n^2=\sum_{k=0}^\infty \left (\sum_{j=1}^n c_{j+k}\right )^2 +
\sum_{k=1}^n \left (\sum_{j=0}^{n-k} c_{j}\right )^2.
\end{equation}

%\begin{example}\label{examp1}
{\bf Example 1.} Let us take at first the case $1/2<\b=1-d<1, \
(0<d<1/2)$. If all $c_k$ have the same sign we know (see
\cite{Giraitis}) that covariances $\g_n$ decay as $n^{-1+2d}$,
$\sum_{n \in \bz}|\g_n|=\infty$ and $A_n^2$ grows as
$n^{3-2\b}=n^{1+2d}$ (long memory in the sense of both definitions).
But if we take $c_k=(-1)^k k^{-\b}, \ k\ge 1, c_0=2$ (such choice of
$c_0$ gives us $\sum_{k=0}^\infty c_k >0$ ), then it is not
difficult to verify that for $n=2m, \ m\ge 1$ all $\g_{2m}$ are
positive, while all $\g_{2m-1}, \ m\ge 1$ are negative but the decay
remains the same: $|\g_n|$ tends to zero as $n^{-1+2d}$, therefore
$\sum_{n \in \bz}|\g_n|=\infty$  and we have long memory in the
sense of Definition \ref{def1}. But  the growth of $A_n^2$ is  only
linear, i.e., as in the case of the short memory in the sense of
Definition \ref{def2}. The most important fact is that in order to
show this we need only the natural condition $\sum_k c_k^2 <\infty$,
simple conditions of alternation $c_k=-c_{k+1}$ and monotonicity
$|c_k|>|c_{k+1}|$, which allow to apply Leibnitz theorem for
convergence of alternating series. Let us denote $C=\left
(\sum_{k=0}^\infty c_k \right )^2$ and apply the particular case of
Toeplitz lemma (see, for example, \cite{Loeve}, p 250 ) which says
that if $b_n=\sum_{k=1}^n a_k \uparrow \infty$ and $x_n\to x$, then
$$
\frac{1}{b_n}\sum_{k=1}^n a_k x_k \to x.
$$
Rewriting the second sum in (\ref{varsn}) as $\sum_{k=0}^{n-1} \left
(\sum_{j=0}^{k} c_{j}\right )^2$ and applying the above formulated
statement with $a_k\equiv 1$ and $x_k=\left (\sum_{j=0}^{k}
c_{j}\right )^2$, we have that
\begin{equation}\label{longmem}
\frac{1}{n}\sum_{k=1}^n  \left (\sum_{j=0}^{n-k} c_{j}\right )^2
 \to C.
\end{equation}
Applying the estimate $|\sum_{j=1}^n c_{j+k}|<|c_{1+k}|$ we easily
get
\begin{equation}\label{longmem1}
\frac{1}{n}\sum_{k=0}^\infty  \left (\sum_{j=1}^{n} c_{j+k}\right
)^2  \to 0.
\end{equation}
From (\ref{varsn}), (\ref{longmem}) and (\ref{longmem1}) we get
$A_n^2 \sim Cn$, which means that the growth of normalizing sequence
is the same as in the case of short memory in the sense of
Definition \ref{def2}. Moreover,  the situation can be even worse.
Taking the same  sequence $c_k=(-1)^k k^{-\b}, \ k\ge 1$, we can
 take $c_0=-\sum_{k=1}^\infty c_k$.  Now  the sequence $A_n^2$ even does
not grow to infinity. Namely, using the condition $\sum_{k=0}^\infty
c_k=0$ and using the property of alternation we can write
\begin{equation}\label{longmem2}
\left | \sum_{j=0}^k c_j \right |=\left | \sum_{j=k+1}^\infty c_j
\right |<|c_{k+1}|.
\end{equation}
Again, rewriting the second sum in (\ref{varsn}) as  earlier  and
applying (\ref{longmem2}), we see that this sum is bounded by the
partial sum of  convergent series $\sum_{k=0}^\infty c_k^2$. Also,
as in (\ref{longmem1}), we see that the first sum in (\ref{varsn})
is bounded by the same convergent series. It remains to show  that
in the case under consideration still we have the relation $\sum_{n
\in \bz}|\g_n|=\infty$. To this aim using the condition
$\sum_{k=0}^\infty c_k=0$  we can write
$$
\g_n=\sum_{k=0}^\infty c_k c_{k+n}=\sum_{k=1}^\infty
c_k(c_{k+n}-c_{n})=I_n^{(1)}+I_n^{(2)},
$$
where
$$
I_n^{(1)}=\sum_{k=1}^\infty c_{2k}(c_{2k+n}-c_{n}), \quad
I_n^{(2)}=\sum_{k=0}^\infty c_{2k+1}(c_{2k+1+n}-c_{n}).
$$
Let us consider the case $n=2m, \ m\ge 1$. Separating the term with
$k=0$ in the second sum and substituting the particular values of
$c_k$ we get
$$
\g_{2m}= \left (\frac{1}{(2m+1)^\b}+\frac{1}{(2m)^\b} \right )
-\frac{1}{(2m)^\b}J_1 +J_2(m)+J_3(m),
$$
where
$$
J_1=\sum_{k=1}^\infty \left (\frac{1}{(2k)^\b}-\frac{1}{(2k+1)^\b}
\right ), \ J_2(m)=\sum_{k=1}^\infty  \frac{1}{(2k(2k+2m))^\b},
$$
and
$$
J_3(m)=\sum_{k=1}^\infty  \frac{1}{((2k+1)(2k+2m+1))^\b}.
$$
The series $J_1$ is alternating, therefore, converging and
$J_1<2^{-\b}<1$, while series in the expressions $J_i(m), \ i=2,3$
are absolutely converging since $2\b>1$. Integral criterion gives us
that both these two series decay as $(2m)^{1-\b}$. Since  for $\b<1$
we have $2\b-1<\b$, therefore $\g_{2m}$ for all $m\ge 1$ are
positive and decay as $(2m)^{1-\b}$, therefore
$\sum_{m=1}^\infty|\g_{2m}|=\infty$. Although we do not need, but it
is possible to show that for $n=2m+1$ all $\g_{2m+1}$ are negative,
have the same order of decay as $\g_{2m}$ and there is monotonicity:
$|\g_{2m+1}|<\g_{2m}$.
%\end{example}

%\begin{example}\label{examp2}
{\bf Example 2.} Now let us consider the case $1<\b=1-d<3/2, \
(-1/2<d<0)$. In this case $\sum_{k=0}^\infty |c_k|<\infty$ and we
have $\sum_{k\in \bz} |\g_k|<\infty$. A linear process with such
filter, according Definition \ref{def1} can be of short memory if
$|\sum_{k=0}^\infty c_k|>0$, or of negative memory if
$\sum_{k=0}^\infty c_k=0$. Again, let us take $|c_k|=k^{-\b}, \ k\ge
1$  with $\b=1-d>1$ and $c_k=0, \ k< 0$. Let us choose
$c_0=-\sum_{k=1}^\infty c_k$, thus we have the case of negative
memory (Definition \ref{def1}). Consider two extreme cases in this
situation. First, let us take all $c_k, \ k\ge 1$ of one sign, let's
say, positive. Using the relation $\sum_{j=0}^{n} c_{j}=-
\sum_{j=n+1}^{\infty} c_{j}$ it is not difficult (we omit the simple
calculations) to show that $A_n^2$ grows as $n^{3-2\b}=n^{1+2d}$ for
$-1/2 <d<0$ ($1<\b<3/2$) and we get that the  limit process for
(\ref{normlinpr2}) is fractional Brownian motion with the  Hurst
parameter $H=1/2+d$. Thus, we have the case of negative memory in
the sense of Definition \ref{def2}, too. If $\b>3/2 \ (d<-1/2)$,
then the sequence $A_n^2$ is bounded and there will be no
convergence of f.d.d.; if $\b=3/2$ than $A_n^2$ will grow
logarithmic, but this growth does not allow to apply Lamperti
theorem (see Theorem 3.4.1 in \cite{Giraitis}).

Now let us consider another extreme case: we take all $c_k$
alternating, that is, $c_k=(-1)^k k^{-\b}, \ k\ge 1$ (with the same
fixed $\b$) and  $c_0=-\sum_{k=1}^\infty c_k$ (negative memory in
the sense of Definition \ref{def1}). In this case it is not
difficult to show  that $A_n^2$ stays bounded.
%In the first case
%(all $c_k, \ k\ge 1$ are of the same sign) the quantity
%$\sum_{j=k+1}^{\infty} c_{j}$ was decaying like $k^{1-\b}$, but for
%$1<\b<3/2$ the sum $\sum_{k=1}^{n} k^{2(1-\b)}$ is growing, while
%now we have (since the series is alternating)
%$$
%|\sum_{j=k}^{\infty} c_{j}|\le |c_k|=k^{-\b}
%$$
%and for $\b>1$ we obtain the convergent series $\sum_k k^{-\b}$.
Thus, we see that under conditions $|c_k|=k^{-1+d}, \ k\ge 1,$ (for
a fixed $-1/2<d<0$) and $c_0=-\sum_{k=1}^\infty c_k$  we can get
that $A_n^2$ grows as $n^{1+2d}$  or stays bounded. It is an
interesting question if it is possible for a fixed $-1/2<d<0$  and
any given $0<\d <1+2d$ to choose the signs of coefficients $c_k$ so
that $A_n^2$ would grow as $n^{\d}$.
%\end{example}

These two examples and considerations before them suggest two
conclusions. First one is that notion of dependence should be
separated from the memory properties, leaving for the expressions
"long-range dependence" and "short-range dependence" only the
meaning that any measure of  dependence is decaying slowly or
quickly, respectively. The second one is that the memory properties
in the case of finite variance should be defined in the same way as
in Definition \ref{def2}, namely,  the case $\a=2$ (only in this
case we must cover two possibilities for a stationary sequence
$\{\xi_i, \ i\in \bz\}$: it can be jointly regularly varying with
the index $2$ or it can be with finite variance) should be included
into Definition \ref{def2}. Such definition allows to treat memory
properties uniquely  in both cases of finite and infinite variances.
Also it is easy to give explanation for such classification. Long
memory means that a stationary process "remember" the past values in
such a way that the volatility of partial sums of this process are
growing more rapidly comparing with the sequence of i.i.d. random
variables, while negative memory means contrary - volatility of
partial sums of this process are growing more slowly. And short
memory means that the partial sums of this process behave in the
same way as  in the case of the sequence of i.i.d. random variables,
which has no memory at all. From this explanation it seems that the
terms "long memory" and "short memory" are not quite logical, if we
would like to leave the term "negative memory", since words "long"
and "short"  has opposite meanings, while from arguments given above
it follows that "long" and "negative" should be as opposite. Thus,
more logical terms would be "positive memory", "zero memory", and
"negative memory", these terms would be coherent with memory
parameters $0<d<1/2, \ d=0, \ -1/2<d<0$ (in the case $\a=2$,  see
\cite{Giraitis}, p. 36; in the case $\a<2$ these intervals would be
$0<d<1-1/\a, \  -1/\a<d<0$). Also these new terms fit well with the
explanation of properties of increments of limit processes for
(\ref{normlinpr2}): this process of partial sums always is with
dependent increments, but in the case of zero memory it "forgets"
this dependence and the limit process is with independent
increments, while in the case of memories (both positive and
negative) the limit process remains with dependent increments.

Considering linear processes with finite variance we saw that for
some filters the normalization constants for partial sums can stay
bounded. The similar situation can be in the case of linear
processes with infinite variance ($0<\a<2$), since typical
normalizing constants are of the form (again we do not take into
account slowly varying functions)
$$
A_n^\a=\sum_{k=0}^\infty \left |\sum_{j=1}^n c_{j+k}\right |^\a +
\sum_{k=1}^n \left |\sum_{j=0}^{n-k} c_{j}\right |^\a,
$$
and similar analysis as in the case $\a=2$ reveals the possibility
for $A_n^\a$ to stay bounded. Therefore it is reasonable to suggest
to call such stationary sequences having strongly negative memory
(memory is so strong that it prevents of growing the volatility of
partial sums of the process). Although for general stationary
sequences with strongly negative memory the problem of limits for
partial sums has no sense (see, for example \cite{Ibragimov}, Ch. 18
), for linear processes even with strongly negative memory this
problem is not trivial if we assume that infinitely many
coefficients of a filter are non-zero.

Defining the memory properties by means of Definition \ref{def2},
the next step will be to clarify what behavior of covariances (the
case of finite variances) and $\a$-covariances (the case of infinite
variances) give us the particular memory property. For general
stationary sequences, without doubt, it is  a difficult problem, but
even for linear processes in the case of finite variance it is not
easy one. If we leave the partition based on convergence or
divergence of the series $\sum_{k\in \bz} |\g_k|$ to separate
long-range dependence and short-range dependence, then one would
guess that short memory (or in new terminology zero memory) will be
under the short-range dependence, namely $\sum_{k\in \bz}
|\g_k|<\infty$ and additional condition $\sum_{k\in \bz} \g_k>0$.
But the cases of positive and negative memories are more
complicated, as Examples 1 and 2 show. For example, from the Example
1 follows that a linear process with long-range dependence can be of
zero memory or even strongly negative memory. This means that the
condition $\sum_{k\in \bz} |\g_k|=\infty$ is insufficient for
positive memory, stronger condition $|\sum_{k\in \bz} \g_k|=\infty$
probably is needed . In the case of infinite variance the situation
is much more complicated, even for linear processes the relation
between memory properties defined in Definition \ref{def2} and the
behavior of $\a$-covariances is not clear.

 At the beginning of this subsection we mentioned
that it is impossible to characterize memory  properties by the
convergence or divergence of the series $ \sum_{k\in \bz} |\rho_n|$,
but one can try the series $ \sum_{k\in \bz} |\rho_n|^{f(\a)}$ with
some function $f(\a)$ for $0<\a\le 2$ with the property $f(2)=1$
 It turns out that such approach is successful with the function $f(\a)=(\a-1)^{-1}$
in the case $1<\a \le 2$ . Namely, from Corollary \ref{corol1} we
have that, in the case $1<\a \le 2$ and $1/\a<\b <1$,
$$
\rho_n\simeq C(\a, \b)n^{1-\b \a}
$$
 therefore,
$ \sum_{k\in \bz} |\rho_n|^{(\a -1)^{-1}}=\infty $ and we have long
memory. If $1<\a \le 2$,  $1<\b $, and $\sum_j c_j\ne 0$, then
$$
\rho_n\simeq C(\a, \b)n^{-\max(\b \a-1, \b)}.
$$
Now $ \sum_{k\in \bz} |\rho_n|^{(\a -1)^{-1}}<\infty $ and,
according results from \cite{Astrauskas}, formulated above, we have
short memory. But these are very particular results, since in
Corollary \ref{corol1} we had investigated the behavior of $\rho_n$
only in the case of regularly varying coefficients of a filter
having constant sign. Further research involving the effect of
alternation is needed, especially the case of negative memory in the
case $\a<2$ remains unclear. In the case $\a=2$ we know that if we
have additional condition that $\sum_j c_j =0$, covariances decay
more quickly comparing with the case $\sum_j c_j \ne 0$ under the
same decay of $|c_j|$. The same effect should be in the case $\a<2$,
but at present we have only conjecture that in the case $1<\a<2, \
\sum_j c_j =0,$ and $c_j=j^{-\b}(1+O(j^{-g(\a)}))$ with some
function $g$, we should get $ \rho_n\simeq C(\a, \b)n^{1-\b \a}.$

\subsection{Memory properties for random  fields}

The memory properties for random fields are less investigated
comparing with the case of processes, even in the case of finite
variance. In this case usually for a stationary random field
$X_{\bk}, \ \bk\in \bz^d$ long memory (which sometimes is used as
synonym for long-range dependence) is defined as the property that
covariance function $\g_\bk := EX_{\B0}X_{\bk}$ is not absolutely
summable: $\sum_{\bk\in \bz^d}|\g_\bk|=\infty$, while summability of
this series means short memory. An alternative approach (albeit not
equivalent) is to define memory properties via spectral density -
 roughly speaking, a random field has long memory if its spectral density
 is unbounded (and has singularity at zero). One of the most popular
 assumptions on the behavior of covariance function is the following
 its growth at infinity
 $$
\g_\bk \sim ||\bk||^{-\b}L(||\bk||)b\left (\frac{\bk}{||\bk||}\right
), \quad 0<\b<d,
 $$
where $L$ is slowly varying at infinity function and $b$ is
non-negative continuous function defined on unit sphere of $\bbr^d$.
Exactly such condition was used in one of the pioneering works on
long range dependence \cite{Dobrushin}, later on it was used with
some modifications (changing the norm, dropping the assumption that
$b$ is non-negative, etc). Similar (in form) condition was used to
describe the growth of spectral density at origin (see, for example,
\cite{Lavansier}). Both such conditions (via covariance function and
spectral density) gives us the so-called isotropic long memory, also
there are papers dealing with non-isotropic long memory of
stationary random fields. But both these two approaches (via
covariance function and spectral density) has the following
shortcomings. As we saw in previous subsections, it is almost
impossible to introduce long and short memories by using substitutes
of covariance such as $\a$-covariance or other similar measures of
dependence in the case of infinite variance. It seems that negative
memory for fields is not introduced even for fields with finite
variance (at least the author have not  seen any paper on this
topic). Therefore it seems quite natural to suggest the same
approach which was suggested for stationary processes and which is
unified both for random processes with finite and infinite variance,
namely, to use the growth of partial sums formed from the random
field under consideration. But before giving the strict definitions
we shall introduce some notations and shall  give some explanations.
For any set $A\subset \bz^d$ let $\#A$ stands for the cardinality of
the set $A$. If for processes we form partial sums by summing the
values of a process over intervals of the increasing length,
situation is more complicated when we pass to random fields, since
now summation is possible over arbitrary sequence of finite
increasing sets $A_n\subset A_{n+1}$ only with requirement that
$\#A_n \to \infty$ or even over some system of set indexed by
multi-indices. If a random field $X_{\bk}, \ \bk\in \bz^d$ consists
of i.i.d. random variables belonging to the domain of the normal
attraction of $\a$-stable law, $0<\a\le 2,$ then the sequence
$b_n^{1/\a}$, where $b_n=\# A_n$, presents the right normalization
for $\sum_{\bk \in A_n}X_{\bk}$ giving $\a$-stable law as a limit
(again, for simplicity of writing we do not take into account slowly
varying functions, since for classification of memory properties
only the exponent in the normalizing sequence is important). Passing
to general stationary random fields we shall require the joint
regular $\a$-variation (or finite variance in the case $\a=2$) and
limit $\a$-stable law (in order to avoid such trivial situation
$X_{\bk}\equiv X$). Then we would like to take the exponent $1/\a$
as characterization of short memory and a boundary value between
long (or, as we wrote, more logical name it would be "positive") and
negative memories, namely, if the normalizing sequence is
$b_n^{1/\a+\d}$ with some $\d>0$ then we have long memory, while if
$\d<0$ then there is negative memory of the field under
consideration. But one can easily notice that such definition of
memory properties would be incorrect in a sense that for a given
random field memory properties may be dependent on the chosen
sequence of sets $A_n$. Looking more carefully at the Definition
\ref{def2}, one can notice that the same situation is for processes:
usually we take summation of the values of a process over intervals
$A_n=\{k \in \bz\: \  1\le k\le n\}$,  but if we take instead of
intervals  the sets $B_n=\{k= 2m \in \bz : \ 1\le m\le n\}$, a
process with short memory (with respect to sets $A_n$) may became of
long memory (with respect to sets $B_n$). Therefore, trying to
define memory for random fields, we must choose some system of sets
in $\bz^d$, and, clearly, in $\bz^d$ there are many possibilities
for such choice, and the classification of memory of stationary
fields, generally speaking, will depend on this choice. Although
from the first glance it seems as unpleasant factor, on the other
hand, such choice gives us more opportunities to investigate memory
properties. It is clear that dependence for random fields is much
more complicated comparing with dependence for processes (it can be
different in different directions). Memory property  is even more
complicated, since, as we noted speaking about processes, the term
"long-range  dependencies" only means that a stationary process (or
a field) "remember" old (or distant in the case of a field) values,
while  memory properties also characterize how a process or a field
remember these values: due to the memory the volatility of partial
sums of the sequence under consideration can be bigger (long memory,
or positive memory in the new terminology) or smaller (negative
memory), comparing with the sequence which has no memory at all
(i.i.d. random variables). Short memory (zero memory) means that the
volatility is the same as in the case of i.i.d. random variables.
Thus, if we suspect that a random field has the so-called isotropic
memory, we will take balls (in Euclidean norm) in $\bz^d$, but, if
we want to look if there is difference in memory properties along
axes, we will take rectangles (or even we can rotate rectangles, if
we suspect that axes with different memory properties do not
coincide with coordinate axes). We shall demonstrate such
possibility by simple example of a linear field, and for the
simplicity of writing we consider the case $d=2$.

%\begin{example}\label{examp3}
{\bf Example 3.} Let us take a linear field (\ref{linfield}) with
innovations having finite variances, and let
\begin{equation}\label{randfield}
Z_{n,m}=\sum_{t=1}^n \sum_{s=1}^m X_{t,s}.
\end{equation}
This means that we take the sets $A_{n,m}=\{(t,s)\in \bz^2: 1\le
t\le n, \ 1\le s \le m \}$ with the cardinality $\#A_{n,m}=nm$ and
we assume that $\min (n,  m) \to \infty$. Let us take the filter of
special structure: $c_{i, j}=a_i b_j$, where $a_i, \ b_i$ are real
and $\sum_{i=0}^\infty a_i^2 <\infty, \ \sum_{j=0}^\infty
b_j^2<\infty $. Although such structure of the filter does not mean
that the random field is factorized into product of two processes,
it turns out that the variance of $Z_{n,m}$ can be factorized, and
this means that memory properties along $t$ and $s$ axis can be
different. One can easily verify that the following formula,
analogous to (\ref{varsn}) is true
\begin{equation}\label{varznm1}
Var Z_{n,m}=(D_1 +D_2)(E_1 +E_2),
\end{equation}
where
 $$
D_1=\sum_{u=0}^\infty A_{n,u}^2, \  D_2=\sum_{u=1}^n A_{n-u}^2, \
E_1=\sum_{v=0}^\infty B_{m,v}^2, \ E_2=\sum_{v=1}^m B_{m-v}^2,
$$
and
$$
A_k=\sum_{t=0}^k a_{t}, \ A_{n,k}=\sum_{t=1}^n a_{t+k}, \quad
B_k=\sum_{t=0}^k b_{t}, \ B_{n,k}=\sum_{t=1}^n b_{t+k}.
$$

Comparing (\ref{varznm1}) with (\ref{varsn}) we see that each factor
in (\ref{varznm1}) has exactly the same structure as the right-hand
side of (\ref{varsn}), only with $a_i$ or $b_i$ instead of $c_i$.
%\end{example}

This example allows us to use the analysis carried for linear
processes and to get that for the random field with such particular
filter we can have all sixteen possible combinations of memory
properties (four for each axis, long (positive), short (zero),
negative, and strongly negative), for example the random field can
have long memory with respect to $t$ (horizontal) axis and negative
memory with respect to $s$ (vertical) axis. To get such combination
it is sufficient to take $a_0=1, \ a_i=i^{-1+d_1}, \ 0<d_1<1/2, \
b_i=i^{-1+d_2}, \ -1/2<d_2<0, i\ge 1, \ b_0= \sum_{i=1}^\infty b_i
$, then the variance of $Z_{n,m}$ will grow as
$n^{1+2d_1}m^{1+2d_2}$ (up to the constant, depending on $d_1,
d_2$). Also, taking alternating coefficients $a_i$ or $b_i$ (or even
both) we can face the situation when the variance will not grow with
respect to one  or even both axes. One more consequence from this
simple example is that in the case where the memory  is
non-isotropic, the cardinality of a set over which is taken partial
summation is not appropriate characteristic: if in the above example
we take $d_1=|d_2|$, then the growth of $Var Z_{n,m}$ will be
proportional to the cardinality of the rectangle of summation (i.e.
$nm$) showing the short memory, while in reality we have long and
negative memories with respect to corresponding coordinate axes.

Based on these considerations we can propose the following general
definition of directional memory, analogous to Definition \ref{def2}
for processes. For simplicity of writing we shall take again the
case $d=2$ (generalization to general case $d>2$ is obvious). Let
$X=(X_{i,j}, i, j \in \bz)$ be a stationary random field with
marginal distribution of $X_{0,0}$  belonging to the domain of
attraction of a stable law with index $0<\a\le 2$, $EX_{0,0}=0$ if
$\a>1$ and $X_{0,0}$ is symmetric if $\a=1$. Let $Z_{n,m}$ be
defined as in (\ref{randfield}).
\begin{definition}\label{def3} We say that a stationary random field
$X$ defined above has directional $(\d_1, \d_2)$-memory, if there
exist slowly varying functions $L_i, \ i=1,2$ such that $A_{n,
m}Z_{n,m}$ converge weakly, as $\min(m,n)\to \infty$, to
non-degenerate bivariate $\a$-stable law, $0<\a \le 2$ and
$$A_{n,m}=\frac{1}{n^{1/\a+\d_1}m^{1/\a+\d_2}L_1(n)L_2(m)}, \quad -\frac{1}{\a}<\d_i <1-\frac{1}{\a}.$$
\end{definition}
Positive value of corresponding $\d$ means positive (long) memory in
corresponding direction, similarly, negative value of $\d$ shows
negative memory, while zero value of $\d$ corresponds to zero
(short) memory. The first step in application of this definition
would  be to prove the limit theorems for linear random fields,
generalizing results in \cite{Astrauskas}.

\section{Extensions, generalizations and open problems}

In the last section we provide several possible extensions or
generalizations of $\a$-covariance function.

\medskip

1) Examining more carefully  the  paper \cite{Paul3} it is possible
to notice  that the restriction of the definition of $\a$-cc
 to $S\a S$ random vectors is superfluous and without difficulty
the notion of $\a$-covariance can be extended from $S\a S$ random
vectors to more general $\a$-stable vectors. In the above cited
paper the reason of this restriction was explained by an example of
$\a$-stable non-symmetric random vector $(X_1, X_2)$ with
independent coordinates, for which both $\a$-covariance and
$\a$-correlation, defined by means of centered random variables
$Y_1$ and $Y_2$ (see construction before formulas (\ref{gapdef}) and
(\ref{acovdef})) generally do not vanish. For a random vector $(X_1,
X_2)$ with independent coordinates and finite second moments
covariance is equal to zero only for centered coordinates, since
$E(X_1-EX_1)(X_2-EX_2)=E(X_1-EX_1)E(X_2-EX_2)$. For $\a$-stable
non-symmetric random vector $(X_1, X_2)$ with independent
coordinates, contrary, $E(Y_1-EY_1)(Y_2-EY_2)=-EY_1EY_2\ne 0$, if
both expectations are non-zero, while $EY_1Y_2=0$, since ${\tilde
\G}$ is concentrated on the axes. Thus, definitions of
$\a$-covariance and $\a$-correlation by formulas (\ref{gapdef}) and
(\ref{acovdef}) (without centering) can be extended to general
$\a$-stable vectors. It is not difficult to verify that all
properties of Proposition \ref{gapprop} remains valid.

\medskip

 2) Since the notion of covariance is defined not only  for Gaussian random
 vectors, but for all vectors having finite second moments (that is, belonging to the domain of normal attraction of a Gaussian law)
 it is natural to try to extend $\a$-covariance for random vectors
 belonging to the domain of normal attraction of  a $S\a S$ random
 vector. We propose to do this in the following way. Let $(\xi_1,
 \xi_2)$ be a random vector satisfying the following condition:
 there exists a finite symmetric measure $\G$ on $S_2$ such that for
 any Borel  set on $S_2$
$$
\lim_{t\to \infty} t^\a P\left (||(\xi_1,
 \xi_2)||>t, \ (\xi_1,
 \xi_2)||(\xi_1,
 \xi_2)||^{-1}\in A \right )=\G (A).
$$
This condition means that the random vector $(\xi_1,
 \xi_2)$ belongs to the domain of normal attraction of  a $S\a S$ random
 vector $\bX=(X_1, X_2)$ with the spectral measure $\G$. We suggest
 to define $\a$-covariance and
 $\a$-correlation of $(\xi_1, \xi_2)$ by means of the measure $\G$, as these quantities
 are defined for $S\a S$ random vector $\bX=(X_1, X_2)$:
\begin{equation}\label{acovardef2}
\rho (\xi_1, \xi_2)=\rho (X_1, X_2)=\int_{S_2}s_1s_2{\G}(ds)
\end{equation}
and similarly for ${\tilde \rho}(\xi_1, \xi_2)$ (taking into account
the first generalization, given above, the assumption of the
symmetry of $\G$ can be dropped). The reason for such definition is
the following. If we consider sums of i.i.d two-dimensional random
vectors with second moment (i.e., $\a=2$), appropriately normalized
by scalars, then the covariance matrix of the limit Gaussian
distribution is the same as that of summands. The similar situation
is in the case $\a<2$, when we consider sums of i.i.d. random
vectors in the domain of $\a$-stable random vector: measure $\G$
which is the main characteristic of summands serves as the spectral
measure of a limit stable distribution, and if we agree that $\G$ is
"responsible" for dependence properties between components of the
limit law, it is natural that the same measure $\G$ defines
dependence for summands. Such extension of the notion of
$\a$-covariance is very useful for linear processes (and fields,
too), since considering linear processes (\ref{linpr}) usually it is
assumed that innovations are only  in  the normal domain of
attraction of $\a$-stable random variable. Thus, let us consider a
linear process
\begin{equation}\label{linpr1}
Z(k)=\sum_{j=0}^\infty c_j\eta_{k-j}, \ k\in \bz,
\end{equation}
where $\{\eta_i, i\in \bz \},$ are i.i.d. random variables belonging
to the normal domain of a standard $S\a S$ random variable with
ch.f. $\exp (-|t|^\a), \ 0<\a < 2$, and a filter $\{c_j, \ j\ge 0\}
$ is such that (\ref{linpr1}) is defined correctly. Then it is easy
to see that finite dimensional distributions of the process $Z$
belong to the normal domain of corresponding distributions of the
process $X$ from (\ref{linpr}), therefore, taking into account
(\ref{acovardef2}), we get
$$
\rho (Z(0), Z(n))= \rho (X(0), X(n)),
$$
and we can use the expressions given in Theorem \ref{thm2}. The same
approach can be used for linear fields, too. Here it is appropriate
to mention that the  codifference for the random vector $(\xi_1,
 \xi_2)$ can be defined directly by the formula (\ref{codifdefgen}),
 but then $\t (\xi_1, \xi_2)$ would not be the same as  $\t (X_1,
 X_2)$. It seems that to calculate $\t (Z(0), Z(n))$ by means of
 (\ref{codifdefgen}) would be rather difficult.
%  Therefore  it is reasonable to suggest
%for the codifference and covariation
% to use the same approach as for $\a$-covariance - for $(\xi_1, \xi_2)$ to define
%these quantities as corresponding quantities for $(X_1, X_2)$.

\medskip

3) Measures of dependence can be considered not only for
finite-dimensional random vectors, but also for random elements with
values in infinite-dimensional Banach (or even more general
topological vector) spaces. Just after appearance of \cite{Paul3}
the author spent a year at Gothenburg university studying infinitely
divisible and stable laws in Banach and Hilbert spaces, the results
of this work were presented in two preprints \cite{Paul11} and
\cite{Paul12}. Part of these results were published later in
\cite{Paul13}, \cite{Paul14}, but part remains unpublished till now.
In \cite{Paul11} (see the end of the paper \cite{Paul13}) the analog
of the correlation matrix $\Lambda_\G$, defined for $k$-dimensional
$S\a S$ random vector (see Proposition 3 in \cite{Paul3})  was
introduced for $S\a S$ random vectors with values in a separable
Banach space. This analog was named pseudo-correlation operator (it
is an operator from $B^*$ (conjugate space of $B$) to $B$, as usual
covariance operator), now, adopting terminology of this paper, it
would be called $\a$-covariance operator. Let us note that with
passing from finite-dimensional space to infinite-dimensional spaces
one faces the principal difficulty: in the case $1<\a$ not all
finite measures on the unit sphere of a Banach space can be spectral
measures of an $\a$-stable measure on $B$, and, as far as I know,
the complete description of such spectral measure in Banach spaces
still is not available. In \cite{Paul11} and \cite{Paul13} some
properties of $\a$-covariance operators of $S\a S$ random vectors
with values in separable Banach spaces, such as compactness and
relations with Gaussian covariance operators,  were described, but
also a lot of open problems were formulated, among them description
of $\a$-stable measures in the space $C(0,1)$ in terms of
$\a$-covariance operators. It is necessary to stress that during the
last decades interest in stationary sequences of random elements in
infinite-dimensional spaces has increased, mainly due to functional
data analysis. Regularly varying time series in Banach spaces are
intensively  investigated, the list of references on this topic is
growing very rapidly, see, for example, \cite{Basrak}, \cite{Hult1},
\cite{Hult2}, \cite{Meinguet} and references therein. We hope that
the notion of $\a$-covariance operator, introduced in \cite{Paul11}
and \cite{Paul13} will be useful in this context, also the approach,
proposed to define memory properties in this paper could be applied
for stationary regularly varying sequences in Banach spaces, only
now one more dimension of complexity will appear - the geometry of
Banach spaces.
 We intend to devote a separate paper to all
these problems.

\medskip

4) As it was mentioned in the introduction, the codifference  can be
defined for general bivariate random vectors, in particular, for
infinitely divisible vectors, containing stable vectors as
particular case. In papers \cite{Rosinski1} and \cite{Rosinski2} it
was demonstrated that the codifference is very useful tool
investigating mixing and ergodicity properties of infinitely
divisible processes. It turns out that it is possible to introduce
the notion of $\a$-covariance for a bivariate infinitely divisible
random vector $\bX=(X_1, X_2)$ without finite variance and with the
L{\' e}vy measure $Q$ (without Gaussian component) in the following
way. If the vector $\bX$ has infinite second moment , the same can
be said about the second moments for the measure $Q$, therefore we
define the analog of $\a$-covariance for $\bX$ (may be it can be
called $Q$-covariance, stressing that dependence between coordinates
of an infinitely divisible vector is reflected by the L{\' e}vy
measure $Q$) as usual covariance for radially re-scaled measure $Q$:
\begin{equation}\label{idcov}
{\kappa} (X_1, X_2)=\int_{\bbr^2 }x_1x_2\frac{Q(dx_1dx_2)}{\max (1,
x_1^2+x_2^2)}.
\end{equation}
Since the L{\' e}vy measure $Q$ has  similar properties as the
spectral measure $\G$ of $\a$-stable measures ($X_1$ and $X_2$ are
independent if and only if (iff) measure $Q$ is concentrated on
axes; coordinates are linearly dependent iff the measure $Q$ is
concentrated on a line going through $0$), such measure of
dependence has the main properties of usual covariance. We note that
in the case of  a stable vector the quantity, defined in
(\ref{idcov}) will be equal to $\a$-covariance defined in
(\ref{acovdef}) multiplied by the constant
$\int_0^\infty(r^{\a-1}\max (1, r^2))^{-1}dr$. Clearly, we can
extend this notion for stationary infinitely divisible sequences and
processes  $X=(X_t),  \ t\in \bz,$ or $t\in \bbr$. If we denote
${\kappa} (t)={\kappa} (X_0, X_t)$, then preliminary considerations
based on Maruyama result from \cite{Maruyama} show that vanishing of
this function as $t\to \infty$ will be necessary condition for
mixing of the process $X$, while (under mild condition on the L{\'
e}vy measure $Q_0$ of $X_0$) vanishing of the codifference function
is necessary and sufficient condition, see Corollary 1 in
\cite{Rosinski1}.

 We provided four possible generalizations or extensions of
the notion of $\a$-covariance and $\a$-covariance function, this
list can be prolonged, but it seems that even the results formulated
above allow to ascertain that these notions can be very useful in
the case of infinite variance and could be good substitute for usual
covariance function. At the same time one must keep in mind that in
case of infinite variance $\a$-covariance is not so universal, as
covariance is in the $L_2$ (Hilbertian) theory. For example , in
linear regression most probably covariation (introduced in
\cite{Kanter} and applied for regression and filtration in \cite{
Miller} and \cite{Cambanis} ) is natural and probably can not be
changed by $\a$-covariance. Also the interesting questions are about
relations of $\a$-covariance  with James orthogonality (the
covariation is directly related, see \cite{Samorod}), association,
mixing, distance covariance. All these questions are subjects for
the future research.

\end{document}